\documentclass{amsart}

\usepackage{xypic}
\input xy
\xyoption{all}
\usepackage{epsfig}
\usepackage{amsthm}
\usepackage{amssymb}
\usepackage{amsmath}
\usepackage{amscd}
\usepackage{color}
\usepackage[T1]{fontenc}

\newcommand{\bg}{\begin{equation}}
\newcommand{\ed}{\end{equation}}
\newcommand{\bga}{\begin{eqnarray}}
\newcommand{\eda}{\end{eqnarray}}
\newcommand{\pf}{\textbf{Proof:\ }}

\def\cbdu{\par{\raggedleft$\Box$\par}}

\newtheorem {Theorem}  {Theorem}

\numberwithin{Theorem}{section}

\newtheorem {Lemma}[Theorem]  {Lemma}

\theoremstyle{definition}
\newtheorem{Definition}[Theorem]{Definition}
\theoremstyle{remark}
\newtheorem{Remark}[Theorem]{\bf Remark}

\newtheorem {Corollary}[Theorem]{\bf Corollary}

\expandafter\chardef\csname pre amssym.def
at\endcsname=\the\catcode`\@ \catcode`\@=11
\def\undefine#1{\let#1\undefined}
\def\newsymbol#1#2#3#4#5{\let\next@\relax
 \ifnum#2=\@ne\let\next@\msafam@\else
 \ifnum#2=\tw@\let\next@\msbfam@\fi\fi
 \mathchardef#1="#3\next@#4#5}
\def\mathhexbox@#1#2#3{\relax
 \ifmmode\mathpalette{}{\m@th\mathchar"#1#2#3}%
 \else\leavevmode\hbox{$\m@th\mathchar"#1#2#3$}\fi}
\def\hexnumber@#1{\ifcase#1 0\or 1\or 2\or 3\or 4\or 5\or 6\or 7\or 8\or
 9\or A\or B\or C\or D\or E\or F\fi}

\font\teneufm=eufm10 \font\seveneufm=eufm7 \font\fiveeufm=eufm5
\newfam\eufmfam
\textfont\eufmfam=\teneufm \scriptfont\eufmfam=\seveneufm
\scriptscriptfont\eufmfam=\fiveeufm

\catcode`\@=\csname pre amssym.def at\endcsname

\newcounter{remark}
\newcommand{\R}{\mathbf{R}}

\renewcommand{\div}{\mbox{div}}

\def \ls {{\lambda_q^{2s}}}

\def  \R   {{\mathbb R}}

\def  \12  {{\frac{1}{2}}}

\def\build#1_#2^#3{\mathrel{\mathop{\kern 0pt#1}\limits_{#2}^{#3}}}

\begin{document}

\title[Regularity propagation for MHD]{Propagation of regularity for the MHD system in optimal Sobolev space}


\author [Mimi Dai]{Mimi Dai}
\address{Department of Mathematics, Stat. and Comp.Sci., University of Illinois Chicago, Chicago, IL 60607,USA}
\email{mdai@uic.edu}





\begin{abstract}

We study the problem of propagation of regularity of solutions to the incompressible viscous non-resistive
magneto-hydrodynamics system. According to scaling, the Sobolev space $H^{\frac n2-1}(\mathbb R^n)\times H^{\frac n2}(\mathbb R^n)$ is critical for the system.
We show that  if a weak solution $(u(t),b(t))$ is in $H^{s}(\mathbb R^n)\times H^{s+1}(\mathbb R^n)$ with $s>\frac n2-1$ at a certain time $t_0$, then it will stay in the space for a short time, provided the initial velocity $u(0)\in H^s(\mathbb R^n)$. In the case that the uniqueness of weak solution in $H^{s}(\mathbb R^n)\times H^{s+1}(\mathbb R^n)$ is known, the assumption of $u(0)\in H^s(\mathbb R^n)$ is not necessary.

\bigskip

KEY WORDS: Magneto-hydrodynamics; propagation of regularity; optimal Sobolev spaces.

\hspace{0.02cm}CLASSIFICATION CODE: 76D03, 35Q35.
\end{abstract}

\maketitle

\section{Introduction}

In this paper we consider the incompressible viscous non-resistive
magneto-hydrodynamics (MHD) system:
\begin{equation}\label{MHD}
\begin{split}
u_t+u\cdot\nabla u-b\cdot\nabla b+\nabla p=&\nu\Delta u,\\
b_t+u\cdot\nabla b-b\cdot\nabla u=&0,\\
\nabla \cdot u=0,  \ \ \nabla\cdot b=&0,
\end{split}
\end{equation}
with the initial conditions
\begin{equation}\label{initial}
u(x,0)=u_0(x),\qquad b(x,0)=b_0(x), \qquad \nabla\cdot u_0=\nabla\cdot b_0=0,
\end{equation}
where $x\in\mathbb{R}^n$ with $n\geq 2$, $t\geq 0$, $u$ is the fluid velocity, $p$ is the fluid pressure, $b$
is the magnetic field, and $\nu>0$ is the viscosity coefficient.
System (\ref{MHD}) describes the dynamics of magnetic field in electrically conducting fluid, for instance, plasmas and salt water. It has been extensively investigated by mathematicians in the last few decades. The quantitative properties of solutions in critical spaces (with respect to the scaling) have arisen great interest. It is known that system (\ref{MHD}) has the following scaling,
\[u_\lambda(x,t)=\lambda u(\lambda x,\lambda^2t), \ \ b_\lambda(x,t)=\lambda b(\lambda x,\lambda^2t), \ \ p_\lambda(x,t)=\lambda^2 p(\lambda x,\lambda^2t)\]
solves the system if $(u(x,t), b(x,t), p(x,t))$ does so, with accordingly scaled initial data. For the Navier-Stokes equation in (\ref{MHD}), $H^{\frac n2-1}(\R^n)$ is scaling invariant (also called being critical). For the magnetic field equation in (\ref{MHD}), one would expect that $H^{\frac n2+1}(\R^n)$ is critical, since it is analogous with the Euler equation. However, the local well-posedness of (\ref{MHD}) in $H^s(\R^n)\times H^s(\R^n)$ for $s>\frac n2$ established by Fefferman, McCormick, Robinson and Rodrigo \cite{FMRR14} suggests that $H^{\frac n2}(\R^n)$, rather than $H^{\frac n2+1}(\R^n)$,  is critical for the magnetic field equation. This ``inconsistency'' with Euler equation may be explained by the fact that the magnetic field equation is linear in $b$, while the Euler equation is nonlinear in $u$. Based on the analysis, one can assert that $H^{\frac n2-1}(\R^n)\times H^{\frac n2}(\R^n)$ is critical for the non-resistive MHD system (\ref{MHD}). 

With insight from the scaling argument, it is natural to seek the optimal space for local well-posedness, which would be $H^{s}(\R^n)\times H^{s+1}(\R^n)$ for $s>\frac n2-1$. Indeed, Fefferman, McCormick, Robinson and Rodrigo \cite{FMRR14} first showed that system (\ref{MHD}) is locally well-posed in $H^{s}(\R^n)\times H^{s}(\R^n)$ for $s>\frac n2$. Later, the same authors \cite{FMRR17} improved the local well-posedness space to  $H^{s+\varepsilon}(\R^n)\times H^{s+1}(\R^n)$ for $s>\frac n2-1$ and a small enough constant $\varepsilon>0$. 
The reason $\varepsilon$ cannot be taken 0 is that the maximal regularity estimate for the Stokes equation from $H^s$ to $L^1(0,T; H^{s+2})$ cannot be obtained. On the other hand, the authors of \cite{CMRR} established local existence for the system in the critical Besov space $B^{\frac n2-1}_{2,1}\times B^{\frac n2}_{2,1}$.

This paper concerns the problem of propagation of regularity of solutions to the non-resistive MHD system (\ref{MHD}) in the optimal Sobolev spaces $H^{s}(\R^n)\times H^{s+1}(\R^n)$ for $s>\frac n2-1$.  The problem of propagation of regularity for the Navier-Stokes equation (NSE) was first studied by Leray in \cite{Leray}. Leray showed that if a weak solution of the NSE is regular at certain time $t_0$, the solution will stay regular for a short time on $(t_0, T)$; and an estimate on $T-t_0$ was obtained. Since the global regularity remains open, such finding regarding regularity propagation is of great interest. Back to the non-resistive MHD system, we will show that if a weak solution is in $H^{s}(\R^n)\times H^{s+1}(\R^n)$ for $s>\frac n2-1$ at some time $t_0>0$, then it will stay in the same space for a short time.
The main ingredient to achieve the goal is the type of maximal regularity estimate for the Stokes equation established in Lemma \ref{le-Sto}. This lemma also has its own interest, for it explains the failure to obtain a solution $u$ of the Stokes equation in $L^1(0,T; H^{s+2})$ provided the initial data $u_0\in H^s(\R^n)$, which is the essential obstacle to remove $\varepsilon$ in the local wellposedness of \cite{CMRR}. The estimate obtained in Lemma \ref{le-Sto} reveals that the obstruction to gain two derivatives in $L^1(0,T)$ is at the initial time. Precisely, the norm of the solution in $L^1(0,T; H^{s+2})$ may blow up on the time interval $(\tau, T)$ like $\log (T/\tau)$ as $\tau\to 0$. 


Our main result states as follows.

\begin{Theorem}\label{thm}
Let $(u, b)$ be a Leray-Hopf weak solution of (\ref{MHD})-(\ref{initial}).  Assume $u_0\in H^s(\R^n)$ with $s>\frac{n}2-1$ and $b_0\in L^2(\R^n)$. In addition, assume that 
$u(t_0)\in H^s(\R^n)$ and $b(t_0)\in H^{s+1}(\R^n)$ for a time $t_0>0$. Then there exists a time \\
$T=T(\|u_0\|_{H^s}, \|u(t_0)\|_{H^s},\|b(t_0)\|_{H^{s+1}})>t_0$ such that 
\[u\in L^\infty(t_0,T; H^s(\R^n))\cap L^2(t_0,T; H^{s+1}(\R^n)),\]
\[b\in L^\infty(t_0,T;H^{s+1}(\R^n)).\]
\end{Theorem}

In the case that the Leray-Hopf weak solution is unique in the space $H^{s}(\R^n)\times H^{s+1}(\R^n)$ with $s>\frac n2-1$, the assumption of $u_0\in H^s(\R^n)$ in Theorem \ref{thm} is not necessary. Namely, we can show the result below.

\begin{Corollary}\label{coro}
Assume $(u_0,b_0)\in L^2(\R^n)\times L^2(\R^n)$. Let $(u, b)$ be the unique Leray-Hopf weak solution of (\ref{MHD})-(\ref{initial}) in $H^{s}(\R^n)\times H^{s+1}(\R^n)$ with $s>\frac n2-1$.   Assume that 
$u(t_0)\in H^s(\R^n)$ and $b(t_0)\in H^{s+1}(\R^n)$ for a time $t_0>0$. Then there exists a time 
$T=T(\|u_0\|_{H^s}, \|u(t_0)\|_{H^s},\|b(t_0)\|_{H^{s+1}})>t_0$ such that 
\[u\in L^\infty(t_0,T; H^s(\R^n))\cap L^2(t_0,T; H^{s+1}(\R^n)),\]
\[b\in L^\infty(t_0,T;H^{s+1}(\R^n)).\]
\end{Corollary}
The justification of Corollary \ref{coro} follows from Theorem \ref{thm} and the uniqueness immediately. Indeed, let $(u,b)$ be the solution of (\ref{MHD})-(\ref{initial}) satisfying assumptions of 
of Corollary \ref{coro}. The basic energy estimate (see Definition \ref{def:weak}) guarantees that
$u\in L^2(0,\infty;H^1(\R^n))$ which implies $u(t)\in H^1(\R^n)$ for almost all $t>0$. Thus, one can pick up a time $\tau_0$ close enough to the initial time 0 such that $u(\tau_0)\in H^1(\R^n)$. If $\frac n2-1<s\leq 1$, $u(\tau_0)\in H^s(\R^n)$ also holds by embedding. Then the uniqueness of weak solution in $H^{s}(\R^n)\times H^{s+1}(\R^n)$ allows us to apply Theorem \ref{thm} by considering $\tau_0$ as the initial time, and the conclusion of propagation of regularity follows right away.

\begin{Remark}
Notice that $H^{\frac n2-1}$ is a critical space for the Navier-Stokes equation. In the high regularity space $H^{s}(\R^n)\times H^{s+1}(\R^n)$ with $s>\frac n2-1$ for the MHD system (\ref{MHD}), there is a good chance that the ``weak-strong'' type of uniqueness holds. In that case, the assumption of $u_0\in H^s(\R^n)$ in Theorem \ref{thm} could be dropped. This issue of uniqueness will be addressed in future work.
\end{Remark}

\bigskip

\section{Preliminaries}
\label{sec:pre}

\subsection{Notation}
\label{sec:notation}
We denote by $A\lesssim B$ an estimate of the form $A\leq C B$ with
some absolute constant $C$, and by $A\sim B$ an estimate of the form $C_1
B\leq A\leq C_2 B$ with some absolute constants $C_1$, $C_2$. We also write
 $\|\cdot\|_p=\|\cdot\|_{L^p}$, and $(\cdot, \cdot)$ stands for the $L^2$-inner product.

\subsection{Littlewood-Paley decomposition}
\label{sec:LPD}
The techniques presented in this paper rely strongly on the frequency localization approach and paradifferential calculus.
Thus we recall the Littlewood-Paley decomposition theory briefly. For a more detailed description on this theory we refer the readers to the books by Bahouri, Chemin and Danchin \cite{BCD} and Grafakos \cite{Gr}. 

Let $\mathcal F$ and $\mathcal F^{-1}$ denote the Fourier transform and inverse Fourier transform, respectively. Define $\lambda_q=2^q$ for integers $q$. A nonnegative radial function $\chi\in C_0^\infty(\R^n)$ is chosen such that 
\begin{equation}\notag
\chi(\xi)=
\begin{cases}
1, \ \ \mbox { for } |\xi|\leq\frac{3}{4}\\
0, \ \ \mbox { for } |\xi|\geq 1.
\end{cases}
\end{equation}
Let 
\bg\notag
\varphi(\xi)=\chi(\frac{\xi}{2})-\chi(\xi)
\ed
and
\begin{equation}\notag
\varphi_q(\xi)=
\begin{cases}
\varphi(\lambda_q^{-1}\xi)  \ \ \ \mbox { for } q\geq 0,\\
\chi(\xi) \ \ \ \mbox { for } q=-1.
\end{cases}
\end{equation}
For a tempered distribution vector field $u$ we define the Littlewood-Paley projection
\begin{equation}\notag
\begin{split}
&h=\mathcal F^{-1}\varphi, \qquad \tilde h=\mathcal F^{-1}\chi,\\
&u_q:=\Delta_qu=\mathcal F^{-1}(\varphi(\lambda_q^{-1}\xi)\mathcal Fu)=\lambda_q^n\int h(\lambda_qy)u(x-y)dy,  \qquad \mbox { for }  q\geq 0,\\
& u_{-1}=\mathcal F^{-1}(\chi(\xi)\mathcal Fu)=\int \tilde h(y)u(x-y)dy.
\end{split}
\end{equation}
By the Littlewood-Paley theory, the following identity
\bg\notag
u=\sum_{q=-1}^\infty u_q
\ed
holds in the distribution sense. Essentially the sequence of the smooth functions $\varphi_q$ forms a dyadic partition of the unit. To simplify the notation, we denote
\bg\notag
u_{\leq Q}=\sum_{q=-1}^Qu_q, \qquad u_{(Q, N]}=\sum_{p=Q+1}^N u_p, \qquad \tilde u_q=\sum_{|p-q|\leq 1}u_p.
\ed

\begin{Definition}
A tempered distribution $u$ belongs to the Besov space $ B_{p, \infty}^{s}$ if and only if
$$
\|u\|_{ B_{p, \infty}^{s}}=\sup_{q\geq-1}\lambda_q^s\|u_q\|_p<\infty.
$$
\end{Definition}
We also note that, 
\[
  \|u\|_{\dot H^s} \sim \left(\sum_{q=-1}^\infty\lambda_q^{2s}\|u_q\|_2^2\right)^{1/2}
\]
for each $u \in  \dot H^s$ and $s\in\R$.

We recall Bernstein's inequality for the dyadic blocks of the Littlewood-Paley decomposition in the following.
\begin{Lemma}\label{le:bern} (See \cite{L}.)
Let $n$ be the space dimension and $r\geq s\geq 1$. Then for all tempered distributions $u$, 
\bg\label{Bern}
\|u_q\|_{r}\lesssim \lambda_q^{n(\frac{1}{s}-\frac{1}{r})}\|u_q\|_{s}.
\ed
\end{Lemma}

\subsection{Bony's paraproduct and commutator}
\label{sec-para}

Bony's paraproduct formula 
\begin{equation}\label{Bony}
\begin{split}
\Delta_q(u\cdot\nabla v)=&\sum_{|q-p|\leq 2}\Delta_q(u_{\leq{p-2}}\cdot\nabla v_p)+
\sum_{|q-p|\leq 2}\Delta_q(u_{p}\cdot\nabla v_{\leq{p-2}})\\
&+\sum_{p\geq q-2} \Delta_q(\tilde u_p \cdot\nabla v_p),
\end{split}
\end{equation}
will be used constantly to decompose the nonlinear terms in energy estimate. 
We will also use the notation of the commutator
\begin{equation} \label{commudef}
[\Delta_q, u_{\leq{p-2}}\cdot\nabla]v_p:=\Delta_q(u_{\leq{p-2}}\cdot\nabla v_p)-u_{\leq{p-2}}\cdot\nabla \Delta_qv_p.
\end{equation}
\begin{Lemma}\label{le-commu}
The commutator satisfies the following estimate, for any $1<r<\infty$
\[\|[\Delta_q,u_{\leq{p-2}}\cdot\nabla] v_p\|_{r}\lesssim \|\nabla u_{\leq p-2}\|_\infty\|v_p\|_{r}.\]
\end{Lemma}

\bigskip

\subsection{Definition of solutions}
\label{sec:sol}
We recall some classical definitions of weak and regular solutions.
\begin{Definition}\label{def:weak}
A weak solution of (\ref{MHD}) on $[0,T]$ (or $[0, \infty)$ if $T=\infty$) is a pair of functions $(u, b)$ in the class 
$$
u, b \in C_w([0,T]; L^2(\mathbb R^3))\cap L^2(0,T; H^1(\mathbb R^3)), 
$$
with $u(0)=u_0, b(0)=b_0$, satisfying (\ref{MHD}) in the distribution sense;
moreover, the following energy inequality 
\begin{equation}\notag
\begin{split}
&\|u(t)\|_2^2+\|b(t)\|_2^2+2\nu\int_{t_0}^t\|\nabla u(s)\|_2^2ds
\leq \|u(t_0)\|_2^2+\|b(t_0)\|_2^2
\end{split}
\end{equation}
is satisfied for almost all $t_0\in(0, T)$ and all $t\in(t_0, T]$.
\end{Definition}

\subsection{Estimate for the Stokes equation}
To provide a general result, we consider the Stokes equation with fractional Laplacian 
\begin{equation}\label{stokes}
u_t+\nu(-\Delta)^\alpha u+\nabla p=f, \ \ \nabla\cdot u=0, \ \ u(0)=u_0
\end{equation}
with $\alpha>0$. We will prove a type of maximal regularity result in $L_t^p$ with delay of an arbitrarily short time.
\begin{Lemma}\label{le-Sto}
Assume $f\in L^r(0,T; H^s(\R^n))$, $u_0\in H^{s}(\R^n)$ and $\nabla \cdot u_0=0$ for $1<r<\infty$ and $s>0$. Let $\varepsilon\in(0,T)$. Then the solution $u$ of the Stokes equation (\ref{stokes}) satisfies 
\begin{equation}\label{est-Sto}
\int_\varepsilon^T\|\nabla^\alpha u(t)\|_{H^{s+\alpha}}\, dt\leq C\nu^{-1}\left(\log{\frac{T}{\varepsilon}}\right) \|u_0\|_{H^s}+C\nu^{-1} T^{1-\frac1r}\|f\|_{L^r(0,T;H^s)}
\end{equation}
for some absolute constant $C$.
\end{Lemma}
\pf The validation of the Littlewood-Paley projection of a function being a test function is discussed in \cite{CD-sqg}. We multiply (\ref{stokes}) by $\Delta_q\left(\lambda_q^{2s+4\alpha}u_q\right)$ and integrate on $\R^n$ to arrive
\begin{equation}\notag
\frac 12\frac {d}{dt} \lambda_q^{2s+4\alpha}\|u_q\|_2^2+ \nu\lambda_q^{2s+4\alpha}\|\nabla^\alpha u_q\|_2^2\leq \lambda_q^{2s+4\alpha}\int_{\R^n}f_q\cdot u_q\, dx.
\end{equation}
Applying H\"older's and Young's inequalities to the flux integral yields
\[\left|\int_{\R^n}f_q\cdot u_q\, dx \right| \leq \|f_q\|_2\|u_q\|_2\leq \frac{\nu}{2}\lambda_q^{2\alpha}\|u_q\|_2^2+\frac 1{\nu \lambda_q^{2\alpha}}\|f_q\|_2^2.\]
It follows from the last two inequalities that
\begin{equation}\notag
\frac {d}{dt} \lambda_q^{2s+4\alpha}\|u_q\|_2^2+ \nu\lambda_q^{2s+4\alpha}\|\nabla^\alpha u_q\|_2^2\leq \frac {1}{\nu}\lambda_q^{2s+2\alpha}\|f_q\|_2^2.
\end{equation}
We then apply Duhamel's principle over $[0,t]$ to get
\begin{equation}\notag
\lambda_q^{2s+4\alpha}\|u_q(t)\|_2^2\leq \lambda_q^{2s+4\alpha}\|u_q(0)\|_2^2e^{-\nu \lambda_q^{2\alpha}t}+\frac {1}{\nu}\lambda_q^{2s+2\alpha}\int_0^t \|f_q(s)\|_2^2e^{-\nu\lambda_q^{2\alpha}(t-s)}\, ds.
\end{equation}
Adding the inequality above for all $q\geq-1$, taking square root and integrating over $(\varepsilon, T]$ gives rise to
\begin{equation}\label{L1-est}
\begin{split}
&\int_\varepsilon^T\left(\sum_{q\geq -1}\lambda_q^{2s+4\alpha}\|u_q(t)\|_2^2\right)^{\frac12}\, dt\\
\leq &\int_\varepsilon^T\left(\sum_{q\geq -1}\lambda_q^{2s+4\alpha}\|u_q(0)\|_2^2e^{-\nu \lambda_q^{2\alpha}t}\right)^{\frac12}\, dt\\
&+\int_\varepsilon^T\left(\sum_{q\geq -1}\frac {1}{\nu}\lambda_q^{2s+2\alpha}\int_0^t \|f_q(s)\|_2^2e^{-\nu\lambda_q^{2\alpha}(t-s)}\, ds\right)^{\frac12}\, dt.
\end{split}
\end{equation}
The estimates of the two terms on the right hand side of (\ref{L1-est}) are shown in the following. We have by some fundamental inequalities that
\begin{equation}\notag
\begin{split}
&\int_\varepsilon^T\left(\sum_{q\geq -1}\lambda_q^{2s+4\alpha}\|u_q(0)\|_2^2e^{-\nu \lambda_q^{2\alpha}t}\right)^{\frac12}\, dt\\
\leq &\int_\varepsilon^T\left(\sum_{q\geq -1}\lambda_q^{2s}\|u_q(0)\|_2^2\right)^{\frac12}\left(\sum_{q\geq -1}\lambda_q^{4\alpha}e^{-\nu \lambda_q^{2\alpha}t}\right)^{\frac12}\, dt\\
\leq &\|u_0\|_{H^s}\int_\varepsilon^T\sum_{q\geq -1}\lambda_q^{2\alpha}e^{-\frac12\nu \lambda_q^{2\alpha}t}\, dt\\
\leq &\|u_0\|_{H^s}\int_\varepsilon^T\left(\sum_{q\geq -1}\lambda_q^{\alpha}e^{-\frac14\nu \lambda_q^{2\alpha}t}\right)^2\, dt.\\ 
\end{split}
\end{equation} 
By Lemma 4.1 from \cite{DQS}, we have 
\[\sum_{q\geq -1}\lambda_q^{\alpha}e^{-\frac14\nu \lambda_q^{2\alpha}t}\lesssim \nu^{-\frac12}t^{-\frac12}.\]
Therefore, we have
\begin{equation}\notag
\int_\varepsilon^T\left(\sum_{q\geq -1}\lambda_q^{2s+4\alpha}\|u_q(0)\|_2^2e^{-\nu \lambda_q^{2\alpha}t}\right)^{\frac12}\, dt
\lesssim \|u_0\|_{H^s}\int_\varepsilon^T\nu^{-1}t^{-1}\, dt\lesssim \nu^{-1}\|u_0\|_{H^s}\log{\frac{T}{\varepsilon}}.
\end{equation} 
To handle the second term on the right hand side of (\ref{L1-est}), we first apply H\"older's inequality,
and then change the order of integration
\begin{equation}\notag
\begin{split}
&\int_\varepsilon^T\left(\sum_{q\geq -1}\frac {1}{\nu}\lambda_q^{2s+2\alpha}\int_0^t \|f_q(\tau)\|_2^2e^{-\nu\lambda_q^{2\alpha}(t-\tau)}\, d\tau\right)^{\frac12}\, dt\\
\leq &(T-\varepsilon)^{\frac12}\left(\int_\varepsilon^T\sum_{q\geq -1}\frac {1}{\nu}\lambda_q^{2s+2\alpha}\int_0^t \|f_q(s)\|_2^2e^{-\nu\lambda_q^{2\alpha}(t-\tau)}\, d\tau\, dt\right)^{\frac12}\\
\leq &(T-\varepsilon)^{\frac12}\left(\int_0^T\int_\tau^T\sum_{q\geq -1}\frac {1}{\nu}\lambda_q^{2s+2\alpha} \|f_q(\tau)\|_2^2e^{-\nu\lambda_q^{2\alpha}(t-\tau)}\, dt\, d\tau\right)^{\frac12}\\
\leq &(T-\varepsilon)^{\frac12}\left(\int_0^T\sum_{q\geq -1}\frac {1}{\nu^2}\lambda_q^{2s} \|f_q(\tau)\|_2^2[1-e^{-\nu\lambda_q^{2\alpha}(T-\tau)}]\, d\tau\right)^{\frac12}\\
\leq &(T-\varepsilon)^{\frac12}\left(\int_0^T\frac {1}{\nu^2} \|f(\tau)\|_{H^s}^2\, d\tau\right)^{\frac12}\\
\lesssim &\frac1\nu T^{1-\frac1{r}}\left(\int_0^T\|f(\tau)\|_{H^s}^r\, d\tau\right)^{\frac1r}
\end{split}
\end{equation}
where $r\geq 2$ is required for the last step. In fact, for $1<r\leq 2$, the estimate can be obtained by duality, see \cite{Krylov}.

Therefore, (\ref{est-Sto}) is obtained by combining the above estimates.
\cbdu

\bigskip

\section{Proof of the main result}
\label{sec-est}


We proceed to prove Theorem \ref{thm} in this section. We will only show a priori estimates satisfied by smooth solutions. A rigorous analysis relies on performing estimates on the Galerkin approximations and then the passage to the limit. 

Formally, multiplying the first equation in (\ref{MHD}) by $\lambda_q^{2s}\Delta_q^2 u$ and the second one by $\lambda_q^{2r}\Delta_q^2 b$, and adding up for all $q\geq -1$ we obtain 
\begin{equation}\label{ineq-uq}
\frac12\frac{d}{dt}\sum_{q\geq -1}\lambda_q^{2s}\|u_q\|_2^2
 +\nu\sum_{q\geq-1}\lambda_q^{2s+2}\|u_q\|_2^2\leq -I_1-I_2,
\end{equation}
\begin{equation}\label{ineq-bq}
\frac12\frac{d}{dt}\sum_{q\geq -1}\lambda_q^{2r}\|b_q\|_2^2 \leq -I_3-I_4,
\end{equation}
with
\begin{equation}\notag
\begin{split}
I_1=&\sum_{q\geq -1}\lambda_q^{2s}\int_{\R^3}\Delta_q(u\cdot\nabla u)\cdot u_q\, dx, \qquad
I_2=-\sum_{q\geq -1}\lambda_q^{2s}\int_{\R^3}\Delta_q(b\cdot\nabla b)\cdot u_q\, dx,\\
I_3=&\sum_{q\geq -1}\lambda_q^{2r}\int_{\R^3}\Delta_q(u\cdot\nabla b)\cdot b_q\, dx,\qquad
I_4=-\sum_{q\geq -1}\lambda_q^{2r}\int_{\R^3}\Delta_q(b\cdot\nabla u)\cdot b_q\, dx,\\
\end{split}
\end{equation}

Using Bony's paraproduct and the commutator notation, $I_1$ is decomposed as
\begin{equation}\notag
\begin{split}
I_1=
&\sum_{q\geq -1}\sum_{|q-p|\leq 2}\lambda_q^{2s}\int_{\R^3}\Delta_q(u_{\leq p-2}\cdot\nabla u_p)u_q\, dx\\
&+\sum_{q\geq -1}\sum_{|q-p|\leq 2}\lambda_q^{2s}\int_{\R^3}\Delta_q(u_{p}\cdot\nabla u_{\leq{p-2}})u_q\, dx\\
&+\sum_{q\geq -1}\sum_{p\geq q-2}\lambda_q^{2s}\int_{\R^3}\Delta_q(u_p\cdot\nabla\tilde u_p)u_q\, dx\\
=&I_{11}+I_{12}+I_{13},
\end{split}
\end{equation}
with 
\begin{equation}\notag
\begin{split}
I_{11}=&\sum_{q\geq -1}\sum_{|q-p|\leq 2}\lambda_q^{2s}\int_{\R^3}[\Delta_q, u_{\leq{p-2}}\cdot\nabla] u_pu_q\, dx\\
&+\sum_{q\geq -1}\sum_{|q-p|\leq 2}\lambda_q^{2s}\int_{\R^3}u_{\leq{q-2}}\cdot\nabla \Delta_q u_p u_q\, dx\\
&+\sum_{q\geq -1}\sum_{|q-p|\leq 2}\lambda_q^{2s}\int_{\R^3}(u_{\leq{p-2}}-u_{\leq{q-2}})\cdot\nabla\Delta_qu_p u_q\, dx\\
=&I_{111}+I_{112}+I_{113}.
\end{split}
\end{equation}
One can see that $I_{112}=0$ due to the fact $\sum_{|p-q|\leq 2}\Delta_qu_p=u_q$ and $\nabla\cdot u_{\leq q-2}=0$.
By the commutator estimate, we obtain
\begin{equation}\notag
\begin{split}
|I_{111}|\leq&\sum_{q\geq -1}\sum_{|p-q|\leq 2}\lambda_q^{2s}\|\nabla u_{\leq p-2}\|_\infty\|u_p\|_2\|u_q\|_2\\
\lesssim &\sum_{q\geq -1}\lambda_q^{2s}\|u_q\|_2^2\sum_{p\leq q}\lambda_p^{\frac n2+1}\|u_p\|_2\\
\lesssim &\sum_{q\geq -1}\lambda_q^{(s+1)\theta}\|u_q\|_2^\theta \lambda_q^{s(2-\theta)}\|u_q\|_2^{2-\theta}\sum_{p\leq q}\lambda_p^{(s+1)\delta}\|u_p\|_2^\delta\lambda_p^{s(1-\delta)}\|u_p\|_2^{1-\delta}\left(\lambda_q^{-\theta}\lambda_p^{\frac n2+1-s-\delta}\right)\\
\lesssim &\sum_{q\geq -1}\lambda_q^{(s+1)\theta}\|u_q\|_2^\theta \lambda_q^{s(2-\theta)}\|u_q\|_2^{2-\theta}\sum_{p\leq q}\lambda_p^{(s+1)\delta}\|u_p\|_2^\delta\lambda_p^{s(1-\delta)}\|u_p\|_2^{1-\delta}\lambda_{p-q}^\theta
\end{split}
\end{equation}
for parameters $\theta$ and $\delta$ satisfying $0<\theta<2$, $0<\delta<1$ and 
\begin{equation}\label{para1}
s\geq \frac n2+1-\theta-\delta.
\end{equation}
It then follows from Young's inequality with $(r_1,r_2,r_3,r_4)\in (1,\infty)^4$ satisfying 
\begin{equation}\label{para2}
\frac1{r_1}+\frac1{r_2}+\frac1{r_3}+\frac1{r_4}=1, \ \ r_1=\frac2\theta, \ \ r_3=\frac2{\delta}
\end{equation}
such that for some $\theta_1>0, \theta_2>0$
\begin{equation}\notag
\begin{split}
|I_{111}|\leq &\frac{\nu}{64}\sum_{q\geq -1}\lambda_q^{2s+2}\|u_q\|_2^2+ C_{\nu}\sum_{q\geq -1}\left(\lambda_q^{2s}\|u_q\|_2^2\right)^{\frac{(2-\theta)r_2}{2}}\\
&+\frac{\nu}{64}\sum_{q\geq -1}\sum_{p\leq q}\lambda_p^{2s+2}\|u_p\|_2^2\lambda_{p-q}^{\theta_1}+C_\nu\sum_{q\geq -1}\sum_{p\leq q}\left(\lambda_p^{2s}\|u_p\|_2^2\right)^{\frac{(1-\delta)r_4}{2}}\lambda_{p-q}^{\theta_2}\\
\leq &\frac{\nu}{32}\sum_{q\geq -1}\lambda_q^{2s+2}\|u_q\|_2^2+C_\nu\left(\sum_{q\geq-1}\lambda_q^{2s}\|u_q\|_2^2\right)^{\frac{(2-\theta)r_2}{2}}
+C_\nu\left(\sum_{q\geq-1}\lambda_q^{2s}\|u_q\|_2^2\right)^{\frac{(1-\delta)r_4}{2}}
\end{split}
\end{equation}
Notice that (\ref{para1}) and (\ref{para2}) imply that $s>\frac n2-1$.

While
\begin{equation}\notag
\begin{split}
|I_{113}|\leq&\sum_{q\geq -1}\sum_{|p-q|\leq 2}\lambda_q^{2s}\|u_{\leq p-2}-u_{\leq q-2}\|_2\|\nabla u_p\|_\infty\|u_q\|_2\\
\lesssim &\sum_{q\geq -1}\lambda_q^{2s+\frac n2+1}\|u_q\|_2^3\\
\lesssim &\sum_{q\geq -1}\lambda_q^{(s+1)\theta}\|u_q\|_2^\theta \lambda_q^{s(3-\theta)}\|u_q\|_2^{3-\theta}\lambda_q^{\frac n2+1-s-\theta}\\
\leq &\frac \nu{32}\sum_{q\geq -1}\lambda_q^{2s+2}\|u_q\|_2^2+C_\nu\left(\sum_{q\geq-1}\lambda_q^{2s}\|u_q\|_2^2\right)^{\frac{3-\theta}{2-\theta}}
\end{split}
\end{equation}
for $s\geq \frac n2+1-\theta$ and $0<\theta<2$. 
Notice that $I_{12}$ and $I_{13}$ can be estimated similarly as $I_{111}$ and $I_{113}$, respectively. Thus 
\begin{equation}\label{est-i1}
I_1\leq \frac{\nu}{8}\|\nabla u\|_{H^s}^2+C_\nu \| u\|_{\dot H^s}^{2+\gamma_1}+C_\nu \| u\|_{\dot H^s}^{2+\gamma_2}
\end{equation}
for $s>\frac n2-1$ and some $\gamma_1, \gamma_2>0$.

Using Bony's paraproduct and the commutator notation, $I_2$ is decomposed as
\begin{equation}\notag
\begin{split}
I_2=
&-\sum_{q\geq -1}\sum_{|q-p|\leq 2}\lambda_q^{2s}\int_{\R^3}\Delta_q(b_{\leq p-2}\cdot\nabla b_p)u_q\, dx\\
&-\sum_{q\geq -1}\sum_{|q-p|\leq 2}\lambda_q^{2s}\int_{\R^3}\Delta_q(b_{p}\cdot\nabla b_{\leq{p-2}})u_q\, dx\\
&-\sum_{q\geq -1}\sum_{p\geq q-2}\lambda_q^{2s}\int_{\R^3}\Delta_q(b_p\cdot\nabla\tilde b_p)u_q\, dx\\
=&I_{21}+I_{22}+I_{23}.
\end{split}
\end{equation}
$I_{21}$ can be estimated as
\begin{equation}\notag
\begin{split}
|I_{21}|\leq & \sum_{q\geq -1}\sum_{|q-p|\leq 2}\lambda_q^{2s+1}\|b_{\leq p-2}\|_\infty\|b_p\|_2\|u_q\|_2\\
\lesssim & \sum_{q\geq -1}\lambda_q^{2s+1}\|b_q\|_2\|u_q\|_2\sum_{p\leq q}\lambda_p^{\frac n2}\|b_p\|_2\\
\lesssim & \sum_{q\geq -1}\lambda_q^{s+1}\|u_q\|_2 \lambda_q^r\|b_q\|_2\sum_{p\leq q}\lambda_p^{r}\|b_p\|_2\lambda_{q-p}^{s-r}\lambda_p^{\frac n2+s-2r}\\
\lesssim & \sum_{q\geq -1}\lambda_q^{s+1}\|u_q\|_2 \lambda_q^r\|b_q\|_2\sum_{p\leq q}\lambda_p^{r}\|b_p\|_2\lambda_{q-p}^{s-r}
\end{split}
\end{equation}
for $\frac n2+s-2r\leq 0$ and $s<r$.
It follows from Young's and Jensen's inequalities that
\begin{equation}\notag
\begin{split}
|I_{21}|\leq &\frac{\nu}{16} \sum_{q\geq -1}\lambda_q^{2s+2}\|u_q\|_2^2+C_\nu\sum_{q\geq -1}\left(\lambda_q^r\|b_q\|_2\sum_{p\leq q}\lambda_p^{r}\|b_p\|_2\lambda_{q-p}^{s-r}\right)^2\\
\leq &\frac{\nu}{16} \sum_{q\geq -1}\lambda_q^{2s+2}\|u_q\|_2^2+C_\nu\left(\sum_{ q\geq-1}\lambda_q^{2r}\|b_q\|_2^2\right)^{1+\gamma}\\
\end{split}
\end{equation}
for some $\gamma>0$. 
We observe
\begin{equation}\notag
\begin{split}
|I_{22}|\lesssim & \sum_{q\geq -1}\lambda_q^{2s}\|b_q\|_2\|u_q\|_2\sum_{p\leq q}\lambda_p^{\frac n2+1}\|b_p\|_2\lesssim |I_{21}|.
\end{split}
\end{equation}
Hence $I_{22}$ shares the same estimate as $I_{21}$.
To handle $I_{23}$, integration by parts followed by H\"older's and Bernstein's inequalities leads to
\begin{equation}\notag
\begin{split}
|I_{23}|=&\left|\sum_{q\geq -1}\sum_{p\geq q-2}\lambda_q^{2s}\int_{\R^3}\Delta_q(b_p\otimes\tilde b_p)\cdot\nabla u_q\, dx\right|\\
\lesssim & \sum_{q\geq -1}\lambda_q^{2s+1}\|u_q\|_2\sum_{p\geq q-4}\|b_p\|_2\|b_p\|_\infty\\
\lesssim & \sum_{q\geq -1}\lambda_q^{2s+1}\|u_q\|_2\sum_{p\geq q-4}\lambda_p^{\frac n2}\|b_p\|_2^2\\
\lesssim & \sum_{q\geq -1}\lambda_q^{s+1}\|u_q\|_2 \sum_{p\geq q-4}\lambda_p^{2r}\|b_p\|_2^2\lambda_{q-p}^{s}\lambda_p^{\frac n2+s-2r}\\
\lesssim & \sum_{q\geq -1}\lambda_q^{s+1}\|u_q\|_2 \sum_{p\geq q-4}\lambda_p^{2r}\|b_p\|_2^2\lambda_{q-p}^{s}
\end{split}
\end{equation}
for $\frac n2+s-2r\leq 0$. It then follows from Young's inequality that
\begin{equation}\notag
\begin{split}
|I_{23}|\leq &\frac{\nu}{16} \sum_{q\geq -1}\lambda_q^{2s+2}\|u_q\|_2^2+C_\nu\sum_{q\geq -1}\left(\sum_{p\geq q-4}\lambda_p^{2r}\|b_p\|_2^2\lambda_{q-p}^{s}\right)^2\\
\leq &\frac{\nu}{16} \sum_{q\geq -1}\lambda_q^{2s+2}\|u_q\|_2^2+C_\nu\left(\sum_{ q\geq-1}\lambda_q^{2r}\|b_q\|_2^2\right)^{2}.
\end{split}
\end{equation}

To conclude, we obtain 
\begin{equation}\label{est-i2}
|I_2|\leq \frac{3\nu}{16}\|\nabla u\|_{H^s}^2+C_\nu \| b\|_{\dot H^r}^{4}+C_\nu \| b\|_{\dot H^s}^{2+\gamma}
\end{equation}
for $\frac n2+s-2r\leq 0$ and $s<r$.

Now we estimate $I_3$ by first decomposing it  as
\begin{equation}\notag
\begin{split}
I_3=
&\sum_{q\geq -1}\sum_{|q-p|\leq 2}\lambda_q^{2r}\int_{\R^3}\Delta_q(u_{\leq p-2}\cdot\nabla b_p)b_q\, dx\\
&+\sum_{q\geq -1}\sum_{|q-p|\leq 2}\lambda_q^{2r}\int_{\R^3}\Delta_q(u_{p}\cdot\nabla b_{\leq{p-2}})b_q\, dx\\
&+\sum_{q\geq -1}\sum_{p\geq q-2}\lambda_q^{2r}\int_{\R^3}\Delta_q(u_p\cdot\nabla\tilde b_p)b_q\, dx\\
=&I_{31}+I_{32}+I_{33},
\end{split}
\end{equation}
with 
\begin{equation}\notag
\begin{split}
I_{31}=&-\sum_{q\geq -1}\sum_{|q-p|\leq 2}\lambda_q^{2r}\int_{\R^3}[\Delta_q, u_{\leq{p-2}}\cdot\nabla] b_pb_q\, dx\\
&-\sum_{q\geq -1}\sum_{|q-p|\leq 2}\lambda_q^{2r}\int_{\R^3}u_{\leq{q-2}}\cdot\nabla \Delta_q b_p b_q\, dx\\
&-\sum_{q\geq -1}\sum_{|q-p|\leq 2}\lambda_q^{2r}\int_{\R^3}(u_{\leq{p-2}}-u_{\leq{q-2}})\cdot\nabla\Delta_qb_p b_q\, dx\\
=&I_{311}+I_{312}+I_{313}.
\end{split}
\end{equation}
One can see that $I_{312}=0$ due to the fact $\sum_{|p-q|\leq 2}\Delta_qb_p=b_q$ and $\nabla\cdot u_{\leq q-2}=0$.
By the commutator estimate and H\"older's inequality, we obtain
 \begin{equation}\notag
\begin{split}
|I_{311}|\leq&\sum_{q\geq -1}\sum_{|p-q|\leq 2}\lambda_q^{2r}\|\nabla u_{\leq p-2}\|_\infty\|b_p\|_2\|b_q\|_2\\
\lesssim &\|\nabla u\|_\infty\sum_{q\geq -1}\lambda_q^{2r}\|b_q\|_2^2\\
\lesssim &\|\nabla u\|_{H^{s+1}}\sum_{q\geq -1}\lambda_q^{2r}\|b_q\|_2^2\\
\end{split}
\end{equation}
since $s>\frac n2-1$. While $I_{313}$ is estimated as
\begin{equation}\notag
\begin{split}
|I_{313}|\leq &\sum_{q\geq -1}\sum_{|p-q|\leq 2}\lambda_q^{2r}\|u_{\leq p-2}-u_{\leq q-2}\|_\infty\|\nabla b_p\|_2\|b_q\|_2\\
\lesssim &\sum_{q\geq -1}\lambda_q^{2r+1}\|u_q\|_\infty\|b_q\|_2^2\\
\lesssim &\|\nabla u\|_\infty \sum_{q\geq -1}\lambda_q^{2r}\|b_q\|_2^2\\
\lesssim &\|\nabla u\|_{H^{s+1}}\sum_{q\geq -1}\lambda_q^{2r}\|b_q\|_2^2.
\end{split}
\end{equation}
Similarly, H\"older's and Bernstein's inequalities applied to $I_{32}$ gives
\begin{equation}\notag
\begin{split}
|I_{32}|=&\left|\sum_{q\geq -1}\sum_{|q-p|\leq 2}\lambda_q^{2r}\int_{\R^3}\Delta_q(u_{p}\cdot\nabla b_{\leq{p-2}})b_q\, dx\right|\\
\leq& \sum_{q\geq -1}\sum_{|q-p|\leq 2}\lambda_q^{2r}\|u_p\|_2\|\nabla b_{\leq p-2}\|_\infty\|b_q\|_2\\
\lesssim & \sum_{q\geq -1}\lambda_q^{2r}\|u_q\|_2\|b_q\|_2\sum_{p\leq q}\lambda_p^{\frac n2+1}\|b_p\|_2\\
\lesssim & \sum_{q\geq -1}\lambda_q^{s+2}\|u_q\|_2\lambda_q^r\|b_q\|_2\sum_{p\leq q}\lambda_p^{r}\|b_p\|_2\lambda_{q-p}^{r-s-2}\lambda_p^{\frac n2-s-1}\\
\lesssim & \sum_{q\geq -1}\lambda_q^{s+2}\|u_q\|_2\lambda_q^r\|b_q\|_2\sum_{p\leq q}\lambda_p^{r}\|b_p\|_2\lambda_{q-p}^{r-s-2}
\end{split}
\end{equation}
since $s>\frac n2-1$. Again it follows from Young's and Jensen's inequalities that
\begin{equation}\notag
\begin{split}
|I_{32}|
\lesssim & \|\nabla u\|_{H^{r+1}} \sum_{q\geq -1}\lambda_q^r\|b_q\|_2\sum_{p\leq q}\lambda_p^{r}\|b_p\|_2\lambda_{q-p}^{r-s-2}\\
\lesssim & \|\nabla u\|_{H^{s+1}} \sum_{q\geq -1}\lambda_q^{2r}\|b_q\|_2^2
\end{split}
\end{equation}
for $r<s+2$.
Integrating by parts for $I_{33}$, we have 
\begin{equation}\notag
\begin{split}
|I_{33}|\leq&\sum_{q\geq -1}\sum_{p\geq q-2}\lambda_q^{2r}\|u_p\|_2\|\tilde b_p\|_2\|\nabla b_q\|_\infty\\
\lesssim &\sum_{q\geq -1}\lambda_q^{2r+\frac n2+1}\|b_q\|_2\sum_{p\geq q}\|u_p\|_2\|b_p\|_2\\
\lesssim &\sum_{q\geq -1}\lambda_q^{r}\|b_q\|_2\sum_{p\geq q}\lambda_p^{s+2}\|u_p\|_2\lambda_p^r\|b_p\|_2\lambda_q^{\frac n2+r+1}\lambda_{p}^{-r-s-2}\\
\lesssim &\sum_{q\geq -1}\lambda_q^{r}\|b_q\|_2\sum_{p\geq q}\lambda_p^{s+2}\|u_p\|_2\lambda_p^r\|b_p\|_2\lambda_q^{\frac n2-s-1}\lambda_{p-q}^{-s-1}\\
\lesssim &\|\nabla u\|_{H^{s+1}}\sum_{q\geq -1}\lambda_q^{r}\|b_q\|_2\sum_{p\geq q}\lambda_p^r\|b_p\|_2\lambda_{p-q}^{-s-1}\\
\end{split}
\end{equation}
since $s> \frac n2-1$. The same routine of applying Young's and Jensen's inequalities gives
\[|I_{33}|\lesssim \|\nabla u\|_{H^{s+1}}\sum_{q\geq -1}\lambda_q^{2r}\|b_q\|_2^2.\]
Thus,  we obtain
for  $s>\frac n2-1$ and $r<s+2$ that
\begin{equation}\label{est-i3}
\begin{split}
|I_{3}|
\lesssim &\|\nabla u\|_{H^{s+1}}\|b_q\|_{\dot H^r}^2.
\end{split}
\end{equation}
Using Bony's paraproduct and the commutator notation, $I_4$ can be written as
\begin{equation}\notag
\begin{split}
I_4=
&-\sum_{q\geq -1}\sum_{|q-p|\leq 2}\lambda_q^{2r}\int_{\R^3}\Delta_q(b_{\leq p-2}\cdot \nabla u_p) b_q\, dx\\
&-\sum_{q\geq -1}\sum_{|q-p|\leq 2}\lambda_q^{2r}\int_{\R^3}\Delta_q(b_{p}\cdot \nabla u_{\leq{p-2}}) b_q\, dx\\
&-\sum_{q\geq -1}\sum_{p\geq q-2}\lambda_q^{2r}\int_{\R^3}\Delta_q(\tilde b_p\cdot \nabla u_p) b_q\, dx\\
=&I_{41}+I_{42}+I_{43}.
\end{split}
\end{equation}
One can observe that $I_{42}$ and $I_{43}$ can be estimated in an analogous way as for $I_{311}$ and $I_{33}$, respectively. Thus we only show the estimate of $I_{41}$, 
\begin{equation}\notag
\begin{split}
|I_{41}|\leq& \sum_{q\geq -1}\sum_{|p-q|\leq 2}\lambda_q^{2r}\|b_{\leq p-2}\|_\infty\|\nabla u_p\|_2\|b_q\|_2\\
\lesssim &\|b\|_\infty \sum_{q\geq -1}\lambda_q^r\|b_q\|_2\lambda_q^{s+1}\|\nabla u_p\|_2\lambda_q^{r-s-1}\\
\lesssim &\|\nabla u\|_{H^{s+1}}\|b\|_\infty \sum_{q\geq -1}\lambda_q^r\|b_q\|_2\\
\lesssim &\|\nabla u\|_{H^{s+1}} \|b_q\|_{\dot H^r}^2
\end{split}
\end{equation}
for $\frac n2<r\leq s+1$. Thus we have for $\frac n2<r\leq s+1$ that
\begin{equation}\label{est-i4}
|I_4|\lesssim \|\nabla u\|_{H^{s+1}} \|b_q\|_{\dot H^r}^2.
\end{equation}

Inequality (\ref{ineq-uq}) along with estimates (\ref{est-i1}) and (\ref{est-i2}) implies that, there exist various constants $C_\nu$ depending on $\nu$ such that
\begin{equation}\label{energy2}
\begin{split}
&\frac{d}{dt}\| u\|_{\dot H^{s}}^2+\nu\|\nabla u\|_{H^s}^2\\
\leq &C_\nu\|u\|_{\dot H^s}^{2+\gamma_1}+C_\nu\|u\|_{\dot H^s}^{2+\gamma_2}
+C_\nu\|b\|_{\dot H^{r}}^{2+\gamma_3}+C_\nu\|b\|_{\dot H^{r}}^4
\end{split}
\end{equation}
with parameters satsifying
\begin{equation}\label{final-cond1}
\frac n2+s-2r\leq 0, \ \quad \frac n2-1<s<r
\end{equation}
and some constants $\gamma_1,\gamma_2,\gamma_3>0$.
Combining estimates (\ref{ineq-bq}), (\ref{est-i3}) and (\ref{est-i4}) gives rise to 
\begin{equation}\notag
\begin{split}
\frac{d}{dt}\|b\|_{\dot H^{r}}^2
\leq &\frac\nu2\|\nabla u\|_{H^s}^2+C_0\|\nabla u\|_{H^{s+1}}\|b\|_{\dot H^{r}}^2
+C_{\nu}\|u\|_{\dot H^s}^{2+\gamma_1}
+C_{\nu}\|b\|_{\dot H^{r}}^{2+\gamma_2}.
\end{split}
\end{equation}
with 
\begin{equation}\label{final-cond2}
\frac n2<r\leq s+1.
\end{equation}
Adding the last two energy inequalities leads to, by dropping similar terms for the sake of simplification
\begin{equation}\label{energy4}
\begin{split}
&\frac{d}{dt}\left(\| u\|_{\dot H^{s}}^2+\|b\|_{\dot H^{r}}^2\right)+\frac\nu2\|\nabla u\|_{H^s}^2\\
\leq &C_{\nu}\left(\| u\|_{\dot H^{s}}^2+\|b\|_{\dot H^{r}}^2\right)^{1+\gamma_1}+C_{\nu}\left(\| u\|_{\dot H^{s}}^2+\|b\|_{\dot H^{r}}^2\right)^{1+\gamma_2}\\
&+C_0\|\nabla u\|_{H^{s+1}}\left(\| u\|_{\dot H^{s}}^2+\|b\|_{\dot H^{r}}^2\right)
\end{split}
\end{equation}
with parameters $r$ and $s$ satisfying (\ref{final-cond1}) and (\ref{final-cond2}). Indeed, we can choose $r=s+1-\delta$ for any $\delta\in[0, \frac12(s-\frac n2+1)]$. For simplicity, we take $r=s+1$ from now on.

Now we pause to estimate $\int_{t_0}^t\|\nabla u(\tau)\|_{H^{s+1}}\, d\tau$ which will appear on the right hand side of (\ref{energy4}) after integration over the time interval $[t_0, t]$. 
First it follows from (\ref{est-Sto}) that for $t\leq t_0+1$ and any $\beta>1$
\begin{equation}\label{L1-1}
\int_{t_0}^t\|\nabla u(\tau)\|_{H^{s+1}}\, d\tau\leq C\nu^{-1} \left(\log{\frac{t}{t_0}}\right) \|u_0\|_{H^{s}}+C\nu^{-1} t^{1-\frac1\beta}\|f\|_{L^\beta(0,t;H^{s})}
\end{equation}
with $f:=-(u\cdot\nabla) u+(b\cdot\nabla)b$. Notice that $s+1>\frac n2$ and hence $H^{s+1}$ is an algebra, we deduce 
\begin{equation}\label{bb}
\|(b\cdot\nabla)b\|_{H^{s}}=\|\nabla\cdot(b\otimes b)\|_{H^{s}}\lesssim \|b\otimes b\|_{H^{s+1}}\lesssim \|b\|_{H^{s+1}}^2\lesssim \|b\|_{H^{s+1}}^2.
\end{equation}
While for the term with $u$,  we have 
\[\|(u\cdot\nabla)u\|_{H^{s}}\lesssim \|u|\nabla|^{s+1}u\|_2\lesssim \|u\|_\infty\|u\|_{H^{s+1}}.\]
Agmon's inequality gives rise to
\[\|u\|_\infty\lesssim \|u\|_{H^{s_1}}^{\theta_1}\|u\|_{H^{\frac n2+\varepsilon_0}}^{1-\theta_1}\]
with $\varepsilon_0=s-\frac n2+1>0$, $\theta_1\in(0,1)$ and $s_1$ satisfying 
\[\frac n2=\theta_1 s_1+(1-\theta_1)(\frac n2+\varepsilon_0).\]
Notice that $s_1=\frac n2+\varepsilon_0(1-\frac1{\theta_1})<\frac n2$. Then by Gagliardo-Nirenberg's interpolation inequality we have 
\[\|u\|_{H^{s_1}}\lesssim \|u\|_{L^2}^{\theta_2}\|u\|_{H^{s+1}}^{1-\theta_2}\]
with $\theta_2\in (0,1)$ satisfying 
\[\frac12=\frac{s_1}{n}+\left(\frac12-\frac{s+1}{n}\right)(1-\theta_2)+\frac{\theta_2}{2}.\]
Putting the last three inequalities together yields
\begin{equation}\label{uu}
\|(u\cdot\nabla)u\|_{H^{s}}\lesssim \|u\|_2^{\theta_1\theta_2}\|u\|_{H^{s+1}}^{2-\theta_1\theta_2}
\end{equation}
with $\theta_1\theta_2=1-\frac{n}{2(s+1)}\in(0,1)$ since $s>\frac n2-1$.
As a consequence of (\ref{bb}) and (\ref{uu}) we have for $\beta=\frac{2}{2-\theta_1\theta_2}>1$ that
\begin{equation}\label{fbs}
\begin{split}
\|f\|_{L^\beta(0,t;H^{s})}^\beta\lesssim & \int_0^t \left(\|b(\tau)\|_{H^{s+1}}^{2}+\|u(\tau)\|_{L^2}^{\theta_1\theta_2}\|u(\tau)\|_{H^{s+1}}^{2-\theta_1\theta_2}\right)^\beta\, d\tau\\
\lesssim & \int_0^t \left(\|b(\tau)\|_{H^{s+1}}^{2\beta}+\|u(\tau)\|_2^{\theta_1\theta_2\beta}\|u(\tau)\|_{H^{s+1}}^{2}\right)\, d\tau.
\end{split}
\end{equation}
Therefore, it follows from (\ref{L1-1}) and (\ref{fbs}) that
\begin{equation}\label{u-L1}
\begin{split}
&\int_{t_0}^t\|\nabla u(\tau)\|_{H^{s+1}}\, d\tau\leq C_\nu \left(\log{\frac{t}{t_0}}\right) \|u_0\|_{H^s}\\
&+C_\nu  (t-{t_0})^{1-\frac1\beta}\left(\int_{t_0}^t \|b(\tau)\|_{H^{s+1}}^{2\beta}+\|u(\tau)\|_2^{\theta_1\theta_2\beta}\|u(\tau)\|_{H^{s+1}}^{2}\, d\tau\right)^{\frac{1}{\beta}}
\end{split}
\end{equation}
with constant $C_\nu$ depending only on $\nu$, $2\beta=\frac{8(s+1)}{2(s+1)+n}$, and $\theta_1\theta_2\beta=\frac{4(s+1)-2n}{2(s+1)+n}$.

In the following, we will proceed a delicate analysis based on (\ref{energy2}),  (\ref{energy4}), (\ref{u-L1}) and a contradiction argument to close the proof of the theorem. We claim that there exists a time $T>t_0$ such that 
\begin{equation}\label{final}
\| u(t)\|_{\dot H^{s}}^2+\|b(t)\|_{\dot H^{s+1}}^2\leq 4(\| u(t_0)\|_{\dot H^{s}}^2+\|b(t_0)\|_{\dot H^{s+1}}^2), \ \ \mbox {for all} \ \ t\in[t_0,T].
\end{equation}
The following notations are adapted:
\begin{equation}\notag
\begin{split}
A(t)= &\| u(t)\|_{\dot H^{s}}^2+\|b(t)\|_{\dot H^{s+1}}^2, \ \ A_0=A(t_0),\\
M_0=&\| u_0\|_2^2+\|b_0\|_2^2,\\
M_1=&C_\nu(4A_0)^{1+\gamma_1}+C_\nu(4A_0)^{1+\gamma_2}+C_\nu(4A_0)^{1+\gamma_3}+C_\nu(4A_0)^2,\\
F(T, A_0,M_0,M_1,\nu,{t_0})=&C_\nu\left(\log\frac{T}{{t_0}}\right) \|u_0\|_{H^s}\\
&+C_\nu (T-{t_0})^{1-\frac{1}{\beta}}\left(A_0^{\beta}T+\nu^{-1}M_0^{\frac12\theta_1\theta_2\beta}(A_0+M_1 T)\right)^{\frac{1}{\beta}}.
\end{split}
\end{equation}
Since $\beta>1$, $F(T, A_0,M_0,M_1,\nu,\varepsilon)$ is increasing in $T$ and $F({t_0},  A_0,M_0,M_1,\nu,{t_0})=0$. Thus, $F$ can be arbitrarily small provided $T$ is arbitrarily close to ${t_0}$.
Indeed, the time $T$ can be chosen as small as that
\begin{equation}\label{choose-T}
e^{F(T,A_0,M_0,M_1,\nu,{t_0})}<2, \ \ \mbox{and} \ \ 2M_1(T-{t_0})/A_0<1.
\end{equation}
 
Take 
\begin{equation}\notag
T_1=\sup\{\tau\in [t_0,T]: A(t)\leq 4A_0 \ \mbox{for all} \ t\in[t_0,\tau] \}.
\end{equation}
Suppose $T_1<T$. Inequality (\ref{energy2}) implies that
\[\int_{t_0}^{T_1}\|\nabla u(t)\|_{H^s}^2\, dt\leq \nu^{-1}(A_0+M_1(T_1-{t_0})).\]
It then follows from (\ref{u-L1}) that
\begin{equation}\notag
\begin{split}
&\int_{t_0}^{T_1}\|\nabla u(\tau)\|_{H^{s+1}}\, d\tau\\
\leq &C_\nu \left(\log{\frac{T_1}{{t_0}}}\right) A_0^{\frac12}
+C_\nu (T_1-{t_0})^{1-\frac1\beta}\left(\int_{t_0}^{T} A_0^{\beta}+M_0^{\frac12\theta_1\theta_2\beta}\|u(\tau)\|_{H^{s+1}}^{2}\, d\tau\right)^{\frac{1}{\beta}}\\
\leq &C_\nu\left(\log{\frac{T_1}{{t_0}}}\right) A_0^{\frac12}
+C_\nu (T_1-{t_0})^{1-\frac1\beta}\left(A_0^{\beta}T_1+\nu^{-1}M_0^{\frac12\theta_1\theta_2\beta}(A_0+M_1T_1)\right)^{\frac{1}{\beta}}\\
=:& F(T_1,A_0,M_0,M_1,\nu,{t_0}).
\end{split}
\end{equation}
Then energy estimate (\ref{energy4}) together with the inequality above implies 
\begin{equation}\notag
\begin{split}
A(T_1)\leq&A({t_0})+M_1(T_1-{t_0})+C_0\int_{t_0}^{T_1}\|\nabla u(\tau)\|_{H^{s+1}}A(\tau)\, d\tau\\
\leq &(A_0+M_1(T_1-{t_0}))\exp\left(\int_{t_0}^{T_1}\|\nabla u(\tau)\|_{H^{s+1}}\, d\tau\right)\\
\leq &(A_0+M_1(T_1-{t_0}))e^{F(T_1,A_0,M_0,M_1,\nu,{t_0})}\\
\leq &2A_0+2M_1(T_1-{t_0})<3A_0,
\end{split}
\end{equation}
where the last two steps follow from the choice of time $T$ as in (\ref{choose-T}).  However, the consequence $A(T_1)<3A_0$ contradicts the definition of $T_1$ and the assumption of $T_1<T$. Therefore $T_1=T$ and (\ref{final}) is justified. In the end, it follows from (\ref{energy4}), (\ref{u-L1}) and (\ref{final}) that
\begin{equation}\notag
\begin{split}
\int_{t_0}^{T}\|\nabla u(t)\|_{H^s}^2\, dt\leq C(\|u({t_0})\|_{H^s}, \|b({t_0})\|_{H^{s+1}}, \nu, {t_0}, n, T),  \\ \int_{t_0}^{T}\|\nabla u(t)\|_{H^{s+1}}\, dt\leq C(\|u({t_0})\|_{H^s}, \|b({t_0})\|_{H^{s+1}}, \nu, {t_0}, n, T),
\end{split}
\end{equation}
for various constants $C$ depending on $\|u(t_0)\|_{H^s}, \|b(t_0)\|_{H^{s+1}}, \nu, {t_0}, n,$ and $T$.
It concludes the proof of Theorem \ref{thm}.

\cbdu

\end{document}